\documentclass[a4paper,3p]{elsarticle}
\usepackage{amsmath,amssymb,mathdots}
\numberwithin{equation}{section}
\newcommand{\ket}[1]{\rvert#1\rangle}

\newmuskip\pFqskip
\mathchardef\pFcomma=\mathcode`, 

\newcommand*\pFq[5]{%
  \begingroup
  \begingroup\lccode`~=`,
    \lowercase{\endgroup\def~}{\pFcomma\mkern\pFqskip}%
  \mathcode`,=\string"8000
  {}_{#1}F_{#2}\biggl[\genfrac..{0pt}{}{#3}{#4};#5\biggr]%
  \endgroup
}

\begin{document}
\begin{frontmatter}
\title{The para-Racah polynomials}
\author{Jean-Michel Lemay\corref{cor1}}
\ead{jean-michel.lemay.1@umontreal.ca}
\author{Luc Vinet\corref{cor2}}
\ead{luc.vinet@umontreal.ca}
\author{Alexei Zhedanov\fnref{label2}\corref{cor1}}
\ead{zhedanov@yahoo.com}
\fntext[label2]{On leave of absence from Donetsk Institute for Physics and Technology, Donetsk 83114, Ukraine}
\address{Centre de recherches math\'ematiques, Universit\'e de Montr\'eal,\\ P.O. Box 6128, Centre-ville Station, Montr\'eal (Qu\'ebec) H3C 3J7}
\date{\today}
\cortext[cor2]{Corresponding author}

\begin{abstract}
New bispectral polynomials orthogonal on a quadratic bi-lattice are obtained from a truncation of Wilson polynomials. Recurrence relation and difference equation are provided. The recurrence coefficients can be encoded in a perturbed persymmetric Jacobi matrix. The orthogonality relation and an explicit expression in terms of hypergeometric functions are also given. Special cases and connections with the para-Krawtchouk polynomials and the dual-Hahn polynomials are also discussed.   
\end{abstract}
\end{frontmatter}
\section{Introduction}

  Hypergeometric orthogonal polynomials have numerous applications. We shall be concerned with polynomials that have $q=1$ as base. The Askey tableau presents a hierarchical organization of these special functions \cite{2010_Koekoek_&Lesky&Swarttouw}. It is comprised of a continuous part and of a discrete one. At the top of the continuous part are the Wilson polynomials expressed in terms of $ _4F_3$ generalized hypergeometric series. A standard truncation condition on the parameters of the Wilson polynomials leads to the Racah polynomials which are orthogonal over a finite set of real points that form a quadratic lattice. The simplest limiting case of the Racah polynomials is that of the Krawtchouk polynomials orthogonal on the linear lattice.

  Two of us have identified a family of orthogonal polynomials that fall outside the Askey scheme \cite{2012_Vinet&Zhedanov_JPhysA_45_265304}. These para-Krawtchouk, as they were called, proved orthogonal over a linear bi-lattice formed by superimposing two linear lattices shifted one with respect to the other. The para-Krawtchouk polynomials naturally arise in quantum transport problems over spin chains \cite{2012_Vinet&Zhedanov_JPhysA_45_265304, Fractional_quantum_spin_chains}.

  We here identify the polynomials that are orthogonal with respect to quadratic bi-lattices. They are obtained from the Wilson polynomials through a novel truncation condition. They have the dual-Hahn polynomials as a special case and the para-Krawtchouk polynomials as a special limit. 

  Consider the Wilson polynomials with parameters $a,b,c,d$ denoted by $W_n(x^2;a,b,c,d)$. They obey the recurrence relation \cite{2010_Koekoek_&Lesky&Swarttouw}
  \begin{align}
  -(a^2+x^2)\tilde{W}_n(x^2)=A_n\tilde{W}_{n+1}(x^2)-(A_n+C_n)\tilde{W}_n(x^2)+C_n\tilde{W}_{n-1}(x^2)
  \end{align}
  where
  \begin{align}
    \tilde{W}_n(x^2)=\frac{W_n(x^2;a,b,c,d)}{(a+b)_n(a+c)_n(a+d)_n}
  \end{align}
  with $(a)_k = a(a+1)\dots(a+k-1)$ the usual Pochhammer symbol and 
  \begin{align} \label{rec_coef_wilson}
  A_n =& \frac{(n+a+b+c+d-1)(n+a+b)(n+a+c)(n+a+d)}{(2n+a+b+c+d-1)(2n+a+b+c+d)} \\
  C_n =& \frac{n(n+b+c-1)(n+b+d-1)(n+c+d-1)}{(2n+a+b+c+d-2)(2n+a+b+c+d-1)}.
  \end{align}
  They also satisfy the difference equation
  \begin{align} \label{dif_eqn_Wilson}
  n(n+a+b+c+d-1)\tilde{W}_n(x^2)=\overline{\mathcal{D}(x)}\tilde{W}_n((x+i)^2)-(\overline{\mathcal{D}(x)}+\mathcal{D}(x))\tilde{W}_n(x^2)+\mathcal{D}(x)\tilde{W}_n((x-i)^2),
  \end{align}
  where $\overline{\mathcal{D}(x)}$ is the complex conjugate of $\mathcal{D}(x)$
  \begin{align}
  \mathcal{D}(x)=\frac{(a+ix)(b+ix)(c+ix)(d+ix)}{(2ix)(2ix+1)}.
  \end{align}
  The Wilson polynomials with parameters $a,b,c,d$, admit an explicit expression given by
  \begin{align}
  \begin{aligned} \label{WilsonPolynomials}
  \tilde{W}_n(x^2;a,b,c,d)=\pFq{4}{3}{-n,n+a+b+c+d-1,a-ix,a+ix}{a+b,a+c,a+d}{1}
		  \equiv \sum_k A_{n,k}\Phi_k(y^2)
  \end{aligned}
  \end{align}
  where
  \begin{align}
  \begin{aligned} \label{ExplicitCoefficients}
  A_{n,k} = \frac{(-n)_k(n+a+b+c+d-1)_k}{(1)_k(a+b)_k(a+c)_k(a+d)_k},\quad \Phi_k(x^2) = (a-ix)_k(a+ix)_k.
  \end{aligned}
  \end{align}
  It is well-known \cite{2010_Koekoek_&Lesky&Swarttouw} that the Wilson polynomials can be reduced to a finite set of $N+1$ orthogonal polynomials if
  \begin{align} \label{truncation_cond}
  A_{N}C_{N+1} = 0.
  \end{align}
  This can be achieved by setting $a+b$, $a+c$, $a+d$, $b+c$, $b+d$ or $c+d$ equal to $-N$. This leads to the Racah polynomials. Another possibility is to take $a+b+c+d-1=-N$, but this introduces a singularity in the denominator of the recurrence coefficients. However, one can get around this problem with the use of limits and obtain new orthogonal polynomials. Our goal here is to study and characterize these polynomials which we shall call para-Racah polynomials. 
  
  This new truncation of Wilson polynomials will be presented in section 2 for odd values of $N$ and the corresponding recurrence relation and difference equation will be obtained. In section 3, we derive an explicit expression for the para-Racah polynomials in terms of hypergeometric functions and compute their weights. The para-Racah polynomials with $N$ even are presented in section 4. In section 5, we discuss special cases and the connections with the para-Krawtchouk polynomials and the dual-Hahn polynomials. 

\section{Recurrence relation and difference equation for $N$ odd}
  
 \subsection{Recurrence relation} 
  Let $N=2j+1$ be an odd integer. The truncation condition \eqref{truncation_cond} with $a+b+c+d-1 = -N$ can be achieved by setting
  \begin{align} \label{TruncationConditions}
  b=-a-j+e_1 t,\quad d=-c-j+e_2 t,
  \end{align}
  and taking the limit $t\to 0$. Here, $a, c$ are free and $e_1, e_2$ are deformation parameters. Inserting \eqref{TruncationConditions} in the recurrence coefficients \eqref{rec_coef_wilson}, it is straightforward to verify that only $A_j$ and $C_{j+1}$ depend on $e_1$ and $e_2$ in the limit $t\rightarrow 0$. Indeed, one finds :
  \begin{align} \label{CoefAn}
  A_n = \begin{cases}
	\dfrac{(n-N)(n+a+c)(n+a-c-j)}{2(2n-N)} &\text{if $n\neq j$,}\\[1em]
	\dfrac{e_1}{e_1+e_2}(j+1)(j+a+c)(a-c) &\text{if $n=j$,}
	\end{cases}
  \end{align}
  \begin{align} \label{CoefCn}
  C_n = \begin{cases}
	\dfrac{n(n-j-1-a+c)(n-N-a-c)(n-j-1)}{(2n-1-N)(2n-N)} &\text{if $n\neq j+1$,}\\[1em]
	\dfrac{e_2}{e_1+e_2}(j+1)(j+a+c)(a-c) &\text{if $n=j+1$.}
	\end{cases}
  \end{align}
  Notice that these new recurrence coefficients are now regular for all $n$. It can be seen from these expressions that the combinations of deformation parameters $e_1, e_2$ that occur are not independent and lead to a single deformation parameter $\alpha$ defined by
  \begin{align} \label{DefAlpha}
  \dfrac{e_1}{e_1+e_2}=\alpha,\quad \dfrac{e_2}{e_1+e_2}=1-\alpha.
  \end{align}
  Hence, the new recurrence coefficients are essentially those of the Wilson polynomials, except for $A_j$ and $C_{j+1}$ which involve the coefficient $\alpha$ :
  \begin{align}
  A_j &= \alpha (j+1)(j+a+c)(a-c) \\
  C_{j+1} &= (1-\alpha)(j+1)(j+a+c)(a-c).
  \end{align}
  It is manifest from \eqref{CoefAn} that the truncation condition $A_NC_{N+1}=0$ is achieved. The resulting recurrence coefficients give rise to a finite set of polynomials  $\tilde{P}_n(x^2)$ that are orthogonal with respect to a discrete measure. Consider the monic version of these polynomials where the leading coefficient is equal to 1 : 
  \begin{align}
  P_n(x^2)= (-1)^n A_1A_2\dots A_{n-1}\tilde{P}_n(x^2) = (x^2)^n+O((x^2)^{n-1}).
  \end{align}
  We shall refer to them as the (normalized) para-Racah polynomials and will denote them by $P_n(x^2;N,a,c,\alpha)$ in general or simply by $P_n(x^2)$ when no confusion can arise. They obey the recurrence relation
  \begin{align} \label{recrel}
  x^2 P_n(x^2) =&\ P_{n+1}(x^2)+(A_n+C_n-a^2)P_n(x^2)+A_{n-1}C_nP_{n-1}(x^2)\\
		\equiv&\ P_{n+1}(x^2)+b_nP_n(x^2)+u_nP_{n-1}(x^2). \notag
  \end{align}
  Using \eqref{CoefAn} and \eqref{CoefCn}, the recurrence coefficients can be expressed as
  \begin{subequations} \label{RecCoef}
  \begin{align}
  b_n =  \begin{cases}
	-\frac{1}{2}[a(a+j)+c(c+j)+n(N-n)] &\text{if $n\neq j$, $j+1$,}\\[1em]
	-a^2-\frac{1}{2}j(1+a-c)(1+a+c+j)+\alpha (a-c)(1+j)(a+c+j) &\text{if $n=j$,} \\[1em]
        -a^2-\frac{1}{2}j(1+a-c)(1+a+c+j)+(1-\alpha)(a-c)(1+j)(a+c+j) &\text{if $n=j+1$,}
	\end{cases}
  \end{align}
  \begin{align}
  u_n =  \begin{cases}
	\frac{n(N+1-n)(N-n+a+c)(n-1+a+c)\left((n-j-1)^2-(a-c)^2\right)}{4(N-2n)(N-2n+2)} &\text{if $n\neq j+1$,}\\[1em]
	\alpha(1-\alpha)(a-c)^2(1+j)^2(a+c+j)^2  &\text{if $n=j+1$.}\\[1em]
	\end{cases} 
  \end{align}
  \end{subequations}
  In order to respect the general positivity condition $u_n > 0$ for $n=1,2,\dots,N$, one must choose one of the two following set of restrictions on $a$, $c$ and $\alpha$ :
  \begin{align}
  \begin{cases}
  -N<a+c<-j+1, \\
  |c-a|>j, \\
  0\leq \alpha \leq 1, 
  \end{cases}\ \text{or}\
  \begin{cases}
  a+c<-N+1\ \text{or}\ 0<a+c, \\
  |c-a|<1, \\
  0\leq \alpha \leq 1. 
  \end{cases}
  \end{align}
  When $\alpha = 1/2$, the coefficients are mirror-symmetric
  \begin{align} \label{PersymmetryCoef}
  &b_n =\ b_{N-n} \qquad    \text{for $n=0,1,\dots,N$}\\
  &u_n =\ u_{N-n+1} \qquad  \text{for $n=1,2,\dots,N$} \notag
  \end{align}
  and for other values $0 \leq \alpha \leq 1$, the middle coefficients $b_j$ and $b_{j+1}$ are perturbed and no longer equal.
  The properties of the para-Racah polynomials can be encoded in the following tri-diagonal or Jacobi matrix with the recurrence coefficients as elements :
  \begin{align} \label{Jmatrix}
  J=
  \begin{pmatrix}
  b_0	&	\sqrt{u}_1	&		&	&
  \\
  \sqrt{u}_1	&	b_1	&	\sqrt{u}_2	&	&
  \\
	  &	\sqrt{u}_2	&	b_2	& \ddots &
  \\
	  & 		& \ddots		& 	\ddots & \sqrt{u}_{N}
  \\
  & & &\sqrt{u}_{N} & b_{N}
  \end{pmatrix}.
  \end{align}
  Its action on the canonical basis $\ket{e_n}$ is 
  \begin{align}
  J\ket{e_n} = \sqrt{u_{n+1}}\ket{e_{n+1}} + b_n\ket{e_n} + \sqrt{u_{n}}\ket{e_{n-1}} 
  \end{align}
  assuming $u_0 = u_{N+1} = 0$. Furthermore, the Jacobi matrix is clearly Hermitian and its eigenvalues define the set of points on which the para-Racah polynomials are orthogonal. One can introduce the eigenbasis
  \begin{align}
  J\ket{s} = x_s\ket{s} 
  \end{align}
  where the eigenvalues are chosen in increasing order $x_0<x_1<\dots<x_N$. It is well known and straightforward to see that the eigenbasis and the canonical one are connected as follows (see e.g. \cite{2012_Vinet&Zhedanov_PhysRevA_85_012323})
  \begin{align}
  \ket{s} = \sum_{n=0}^{N} \dfrac{\sqrt{w_s}\ P_n(x_s)}{\sqrt{u_1 \dots u_n}}\ket{e_n}
  \end{align}
  where $P_n$ denotes the para-Racah polynomials and $w_s$ their weights. The mirror-symmetry \eqref{PersymmetryCoef} implies that the matrix J is persymmetric which means that it is symmetric with respect to the main anti-diagonal:
 \begin{align} \label{persymmetry2}
  RJR=J
  \end{align}
  with
  \begin{align}
  R = 
  \begin{pmatrix}
    & & & 1  \\
    & & 1&   \\
    & \iddots & &  \\
    1& & &  \\
  \end{pmatrix}.
  \end{align}
  This persymmetry property has been widely studied especially in the context of inverse spectral problems \cite{2004_Gladwell, DEBOOR1978245, 1983_Golub&VanLoan_MatrixComp, 0266-5611-3-4-010}. 
  
 \subsection{Difference equation}
  By inserting the parametrization \eqref{TruncationConditions} in the difference equation of the Wilson polynomials \eqref{dif_eqn_Wilson}, it is easily seen that the limit $t\to0$ is trivial since there are no parameters in the denominator. Hence, the para-Racah polynomials obey the same difference equation as the Wilson polynomials :
  \begin{align}
  -n(N-n)P_n(x^2)=\overline{\mathcal{D}(x)}P_n((x+i)^2)-(\overline{\mathcal{D}(x)}+\mathcal{D}(x))P_n(x^2)+\mathcal{D}(x)P_n((x-i)^2)
  \end{align}
  with
  \begin{align}
  \mathcal{D}(x)=\frac{(a+ix)(-a-j+ix)(c+ix)(-c-j+ix)}{(2ix)(2ix+1)}. 
  \end{align}
  Remark that the spectrum is doubly-degenerate since $P_n$ and $P_{N-n}$ are eigenfunctions with the same eigenvalue $-n(N-n)$. These polynomials are hence bispectral but they do not belong to classical orthogonal polynomials in the usual sense.

\section{Explicit expression and orthogonality relation for $N$ odd}

\subsection{Explicit expression}
  The explicit expression for the para-Racah polynomials can be found via a limit of the Wilson polynomials. Using the parametrization \eqref{TruncationConditions} in the coefficients \eqref{ExplicitCoefficients} yields
  \begin{align} \label{ANK}
  A_{n,k} = \frac{(-n)_k(n-N+(e_1+e_2)t)_k}{(1)_k(-j+e_1t)_k(a+c)_k(a-c-j+e_2t)_k}.
  \end{align}
  The monic para-Racah polynomials can be defined as 
  \begin{align} \label{defPR}
  P_{n}(x^2)=\eta_n\sum_k \lim_{t\to0} A_{n,k}\Phi_{k}(x^2) 
  \end{align}
  where $\eta_n$ simply is a normalization factor to ensure monicity. Using \eqref{DefAlpha} and \eqref{ANK}, a simple calculation gives  
  \begin{align}
  \lim_{t\to0} A_{n,k}=\begin{cases}
  \dfrac{(-n)_k(n-N)_k}{(1)_k(-j)_k(a+c)_k(a-j-c)_k}                              &\text{if $k\le j$ and $k\le n$},\\[1.1em]
  \dfrac{\alpha^{-1}(-n)_k(n-N)_{N-n}(1)_{k-1+n-N}}{(1)_k(-j)_j(1)_{k-j-1}(a+c)_k(a-j-c)_k} &\text{if $k>j$ and $k\le n$},\\[1.1em]
  \hfill 0\hfill &\text{otherwise}.
  \end{cases}
  \end{align}
  For $n\leq j$, the sum \eqref{defPR} corresponds to the hypergeometric function
  \begin{align}
  P_n(x^2)=\eta_n&\pFq{4}{3}{-n,n-N,a-ix,a+ix}{-j,a+c,a-c-j}{1},
  \end{align}
  whereas for $n>j$, the sum splits in two hypergeometric functions for $k\leq j$ and $k>j$ respectively :
  \begin{align}
  \begin{aligned}
  P_n(x^2)=\eta_n&\pFq{4}{3}{-n,n-N,a-ix,a+ix}{-j,a+c,a-c-j}{1}\\ &+ \eta_n\frac{(-n)_{j+1}(n-N)_{N-n}(a-ix)_{j+1}(a+ix)_{j+1}(1)_{n-j-1}}{\alpha(1)_{j+1}(-j)_j(a+c)_{j+1}(a-c-j)_{j+1}}\\ &\ \times\pFq{4}{3}{-n+j+1,n-j,a+j+1-ix,a+j+1+ix}{j+2,a+c+j+1,a-c+1}{1}.
  \end{aligned}
  \end{align}
  (Take note that even though there is a negative integer $-j$ in the bottom parameters of the hypergeometric functions, the truncation of the series occurs with $-n$ for $n\le j$ and with $n-N$ for $n>j$ before a zero appears in the denominator). The normalization factor is given by
  \begin{align}
  \eta_n= \begin{cases}
  \dfrac{(1)_n(-j)_n(a+c)_n(a-c-j)_n}{(-n)_n(n-N)_n} &\text{if $n\leq j$},\\[1.1em]
  \dfrac{\alpha(1)_n(-j)_j(1)_{n-j-1}(a+c)_n(a-c-j)_n}{(-n)_n(n-N)_{N-n}(1)_{2n-1-N}} &\text{if $n>j$}.
  \end{cases}
  \end{align}
  
\subsection{Orthogonality relation}  
  For a finite set of $(N+1)$ orthogonal polynomials of degree $n=0,\dots,N$, the zeros of the characteristic polynomial $P_{N+1}$ of degree $N+1$ define the grid on which the polynomials are orthogonal. Here, it can be computed via 
  \begin{align} \label{CharacteristicQ}
  P_{N+1}(x^2)=\eta_{N+1}\sum_k\lim_{t\to0} t A_{N+1,k}\Phi_{k}(x^2),
  \end{align}
  where an additionnal $t$ is necessary to obtain non-zero coefficients. The sum can be carried with the use of the Saalsch\"utz summation formula \cite{2010_Koekoek_&Lesky&Swarttouw} to get 
  \begin{align}
  P_{N+1}(x^2) = \prod_{s=0}^j \left((a+s)^2+x^2\right)\left((c+s)^2+x^2\right).
  \end{align}
  Hence, the para-Racah polynomials are orthogonal on the quadratic bi-lattice
  \begin{align} \label{bi-lattice}
  \begin{aligned}
    x_{2s}= -(s+a)^2, \quad s=0,\dots,j, \\
    x_{2s+1}= -(s+c)^2, \quad s=0,\dots,j.
  \end{aligned}
  \end{align} 
  From the standard theory of orthogonal polynomials \cite{1978_Chihara}, the discrete weights can be obtained via the formula 
  \begin{align} \label{firstweight}
  w_s = \dfrac{u_1 \dots u_N}{P_N(x_s)P'_{N+1}(x_s)}, \quad s=0,1,\dots,N 
  \end{align}
  and the orthogonality relation is   
  \begin{align} \label{OR}
  \sum_{s=0}^N P_n(x_s)P_m(x_s)w_s = u_1 \dots u_n \delta_{nm}. 
  \end{align}
  There is however a simpler procedure to compute the weights that has been explained in \cite{WIP_Genest&Vinet&Zhedanov_isospectral_perturbations}. Recall that the persymmetry property \eqref{PersymmetryCoef} is observed when $\alpha=1/2$. It thus follows that in this case, the polynomial $P_N(x^2)$ takes the following values at the spectral points \cite{2012_Vinet&Zhedanov_PhysRevA_85_012323, DEBOOR1978245} :
  \begin{align} 
   P_N(x_s)=\sqrt{u_1 \dots u_N} (-1)^{s+1}.
  \end{align}
  formula \eqref{firstweight} then reduces to
  \begin{align}
   \tilde{w}_s = \frac{\sqrt{u_1 \dots u_N}}{|P'_{N+1}(x_s)|}. 
  \end{align}
  The corresponding positive weights can straightforwardly be computed
  \begin{align} \label{persymmetricweights}
  \begin{aligned}
  \tilde{w}_{2s}&=\frac{\kappa_N (-j)_s(2a)_s(a+1)_s(a-c-j)_s(a+c)_s}{(a+c)_{j+1}(c-a)_{j+1}(2a+1)_j s!(a)_s(2a+1+j)_s(a-c+1)_s(a+c+j+1)_s}, \\
  \tilde{w}_{2s+1}&=\frac{-\kappa_N (-j)_s(2c)_s(c+1)_s(c-a-j)_s(a+c)_s}{(a+c)_{j+1}(a-c)_{j+1}(2c+1)_j s!(c)_s(2c+1+j)_s(c-a+1)_s(a+c+j+1)_s}
  \end{aligned}
  \end{align}
  with
  \begin{align}
  \kappa_N= \frac{(a-c-j)_N(a+c)_N}{2(-1)^{j+1}\binom{2j}{j}j!}.
  \end{align}
  It has also been shown in \cite{WIP_Genest&Vinet&Zhedanov_isospectral_perturbations} that the weights for a general $\alpha$ are related to those of the mirror-symmetric case by a simple multiplicative factor 
  \begin{align}
   w_s = \emph{const}(1+\beta (-1)^s)\tilde{w_s}
  \end{align}
  where $\beta$ is a real parameter. It is a simple matter to identify this parameter by comparing \eqref{firstweight} and \eqref{persymmetricweights} for a fixed value of $N$. One easily verifies that $\beta = 1-2\alpha$. The general weights are then given by
  \begin{align} \label{odd_weights}
  \begin{aligned}
  w_{2s} &= 2(1-\alpha)\tilde{w}_{2s}, \\ 
  w_{2s+1} &= 2\alpha \tilde{w}_{2s+1}. 
  \end{aligned}
  \end{align}
  The weights have the special property that
  \begin{align} \label{weight_property}
  \sum_{s=0}^{j} w_{2s} = 1-\alpha, \qquad \sum_{s=0}^{j} w_{2s+1} = \alpha,
  \end{align}
  which generalizes a known result when $\alpha = 1/2$ for persymmetric Jacobi matrices \cite{WIP_Genest&Vinet&Zhedanov_isospectral_perturbations}.

\section{Even $N$ case}
  We have so far only considered the truncation of the Wilson polynomials for odd $N$. The even case is treated analogously. We shall thus only summarize how this is done and provide the results when $N=2j$. First of all, the truncation condition \eqref{truncation_cond} is now achieved by the parametrization
  \begin{align} \label{TruncationConditions2}
    b=-a-j+e_1 t,\quad d=1-c-j+e_2 t,
  \end{align}
  where, as before, the limit $t\to 0$ is to be taken. With the help of \eqref{TruncationConditions2} and \eqref{DefAlpha}, the recurrence coefficients \eqref{rec_coef_wilson} become
  \begin{align} \label{CoefAn2}
    A_n = \begin{cases}
	  \dfrac{(n-N)(n+a+c)(n+a-c-j+1)}{2(2n+1-N)} &\text{if $n\neq j$,}\\[1em]
	  \alpha j(j+a+c)(c-a-1) &\text{if $n=j$,}
	  \end{cases}
    \end{align}
    \begin{align} \label{CoefCn2}
    C_n = \begin{cases}
	  \dfrac{n(n-N-a-c)(n-j-1+c-a)}{2(2n-1-N)} &\text{if $n\neq j$,}\\[1em]
	  (1-\alpha)j(j+a+c)(c-a-1) &\text{if $n=j$.}
	  \end{cases}
    \end{align}
  The recurrence relation of the para-Racah polynomials is \eqref{recrel} with
    \begin{subequations} \label{RecCoef2}
    \begin{align}
    b_n = \frac{(n-N)(n+a+c)(n+a-c-j+1)}{2(2n+1-N)}+\frac{n(n-N-a-c)(n-j-1+c-a)}{2(2n-1-N)}-a^2 
    \end{align}
  and
    \begin{align}
    u_n =  \begin{cases}
	  \frac{n (2 j-n+1) (a+c+n-1) (a-c+j-n+1) (-a+c+j-n) (a+c+2 j-n)}{4 (N-2 n+1)^2} &\text{if $n\neq j,\ j+1$,}\\[1em]
	  -\frac{1}{2} (1-\alpha) j (j+1) (a-c) (a-c+1) (a+c+j-1) (a+c+j) &\text{if $n=j$,}\\[1em]
	  -\frac{1}{2} \alpha  j (j+1) (a-c) (a-c+1) (a+c+j-1) (a+c+j) &\text{if $n=j+1$.}
	  \end{cases} 
    \end{align}
    \end{subequations}
  The positivity condition $u_n>0$ for $n=1,2,\dots,N$ now requires one of the following set of constraints on $a$, $c$ and $\alpha$ :
  \begin{align}
  \begin{cases}
  -j<a+c<-j+1,\\
  a-c<-j\ \text{or}\ c-a<-j+1,\\
  0\le\alpha\le1, 
  \end{cases}\ \text{or}\
  \begin{cases}
  a+c<-N+1\ \text{or}\ 0<a+c,\\
  0<c-a<1,\\
  0\le\alpha\le1.   
  \end{cases}
  \end{align}
  The coefficients are again mirror-symmetric for $\alpha = 1/2$, but in this case, it is $u_j$ and $u_{j+1}$ that are perturbed for other values of $\alpha$ with $0\le \alpha \le 1$.
 
  Again, the difference equation for the Wilson polynomials does not change under the parametrization \eqref{TruncationConditions2} and can be written as 
  \begin{align}
  -n(N-n)P_n(x^2)=\overline{\mathcal{D}(x)}P_n((x+i)^2)-(\overline{\mathcal{D}(x)}+\mathcal{D}(x))P_n(x^2)+\mathcal{D}(x)P_n((x-i)^2)
  \end{align}
  with
  \begin{align}
  \mathcal{D}(x)=\frac{(a+ix)(-a-j+ix)(c+ix)(1-c-j+ix)}{(2ix)(2ix+1)}. 
  \end{align}
  The spectrum of the corresponding Jacobi matrix is the same for $N$ odd and $N$ even and each eigenvalue is doubly-degenerate, except for the level $n=j$ which is non-degenerate when $N=2j$. The para-Racah polynomials can be expressed as a sum of the form \eqref{defPR} with 
    \begin{align}
    \lim_{t\to0} A_{n,k}=\begin{cases}
    \dfrac{(-n)_k(n-N)_k}{(1)_k(-j)_k(a+c)_k(a-j-c+1)_k}                              &\text{if $k\le j$ and $k\le n$},\\[1.1em]
    \dfrac{\alpha^{-1}(-n)_k(n-N)_{N-n}(1)_{k-1+n-N}}{(1)_k(-j)_j(1)_{k-j-1}(a+c)_k(a-j-c+1)_k} &\text{if $k>j$ and $k\le n$},\\[1.1em]
    \hfill 0\hfill &\text{otherwise.}
    \end{cases}
    \end{align}
  For $n\leq j$, the sum corresponds to the hypergeometric function
    \begin{align}
    P_n(x^2)=\eta_n&\pFq{4}{3}{-n,n-N,a-ix,a+ix}{-j,a+c,a-c-j+1}{1},
    \end{align}
    whereas for $n>j$, the sum splits in two hypergeometric functions for $k\leq j$ and $k>j$ respectively :
    \begin{align}
    \begin{aligned}
    P_n(x^2)=\eta_n&\pFq{4}{3}{-n,n-N,a-ix,a+ix}{-j,a+c,a-c-j+1}{1}\\ &+ \eta_n\frac{(-n)_{j+1}(n-N)_{N-n}(a-ix)_{j+1}(a+ix)_{j+1}(1)_{n-j}}{\alpha(1)_{j+1}(-j)_j(a+c)_{j+1}(a-c-j+1)_{j+1}}\\ &\ \times\pFq{4}{3}{-n+j+1,n-j+1,a+j+1-ix,a+j+1+ix}{j+2,a+c+j+1,a-c+2}{1}.
    \end{aligned}
    \end{align}
    (As for the $N$ odd case, the truncation of the hypergeometric series occurs before the negative integer $-j$ in the bottom parameters produces a zero). The normalization factor is given by
    \begin{align}
    \eta_n= \begin{cases}
    \dfrac{(1)_n(-j)_n(a+c)_n(a-c-j+1)_n}{(-n)_n(n-N)_n} &\text{if $n\leq j$},\\[1.1em]
    \dfrac{\alpha(1)_n(-j)_j(1)_{n-j-1}(a+c)_n(a-c-j+1)_n}{(-n)_n(n-N)_{N-n}(1)_{2n-1-N}} &\text{if $n>j$}.
    \end{cases}
    \end{align}
  Finally, the orthogonality relation for $N=2j$ is again \eqref{OR} on the bi-lattice
  \begin{align} \label{bi-lattice2}
  \begin{aligned}
    x_{2s}&= -(s+a)^2, \quad s=0,\dots,j, \\
    x_{2s+1}&= -(s+c)^2, \quad s=0,\dots,j-1,
  \end{aligned}
  \end{align}   
  with the weights now given by 
  \begin{align} \label{evenw}
  \begin{aligned}
  w_{2s}  &=\frac{ 2(1-\alpha)\kappa_N (-j)_s(2a)_s(a+1)_s(a-c-j+1)_s(a+c)_s}{j! (a+c)_j(c-a)_j(2a+1)_j(a)_s s! (2a+1+j)_s(a-c+1)_s(a+c+j)_s}, \\ 
  w_{2s+1}&=\frac{-2 \alpha   \kappa_N (-j+1)_s(2c)_s(c+1)_s(c-a-j)_s(a+c)_s}{(j-1)! (a+c)_{j+1}(a-c)_{j+1}(2c+1)_{j-1}(c)_s s! (2c+j)_s(c-a+1)_s(a+c+j+1)_s}, 
  \end{aligned}
  \end{align}
  with
  \begin{align}
  \kappa_N = \frac{(a-c-j+1)_N(a+c)_N}{(-1)^{j}\binom{2j}{j}}.
  \end{align}
  The above weights satisfy 
  \begin{align}
   \sum_{s=0}^{j} w_{2s} = 1-\alpha, \qquad \sum_{s=0}^{j-1} w_{2s+1} = \alpha,
  \end{align}
  that compares to \eqref{weight_property} with a different range in the second sum. 
  
\section{Special cases}

  \subsection{Connection with the para-Krawtchouk polynomials}
  The para-Racah polynomials that we have introduced here are orthogonal on a quadratic bi-lattice given by \eqref{bi-lattice} for $N$ odd and \eqref{bi-lattice2} for $N$ even. It is possible to further deform these bi-lattices into a linear bi-lattice. Take $a$ and $c$ as follows
  \begin{align} \label{reparac}
    a(\theta) = \frac{\theta - \Delta}{2}, \qquad c(\theta) = \frac{\theta + \Delta}{2},
  \end{align}
  it is straightforward to verify that 
  \begin{align}
  \begin{aligned}
  \lim_{\theta\to\infty} \frac{-2(x_{2s}+a(\theta)^2)}{\theta}   = 2s, \quad
  \lim_{\theta\to\infty} \frac{-2(x_{2s+1}+a(\theta)^2)}{\theta} = 2s + 2\Delta, 
  \end{aligned}
  \end{align}
  where $x_i$ is given by \eqref{bi-lattice} or \eqref{bi-lattice2}. This linear bi-lattice corresponds to the orthogonality set of the para-Krawtchouk polynomials \cite{2012_Vinet&Zhedanov_JPhysA_45_265304} with parameter $2\Delta$. The same procedure can be applied to \eqref{recrel} to obtain the recurrence relation of the para-Krawtchouk polynomials. Indeed, consider the polynomials $Q_n(y)$ related to the para-Racah polynomials $P_n(x^2; N,a,c,\alpha)$ by 
  \begin{align}
  x^2 = -\frac{\theta}{2}y- \left(\frac{\theta-\Delta}{2}\right)^2, \quad Q_n(y) = \left(-\frac{\theta}{2}\right)^{-n}P_n\left(x^2; N,\frac{\theta-\Delta}{2},\frac{\theta+\Delta}{2},\frac{1}{2}\right).
  \end{align}
  Inserting this in the recurrence relation \eqref{recrel} yields
  \begin{align} \label{LimitingRel}
  y Q_n(y) = Q_{n+1}(y) -\frac{2b_n+2a(\theta)^2}{\theta}\ Q_n(y) + \frac{4 u_n}{\theta^2}\ Q_{n-1}(y). 
  \end{align}
  Taking the limit $\theta\to\infty$ with the recurrence coefficients $b_n$ and $u_n$ given by \eqref{RecCoef} gives for $N$ odd
  \begin{align}
  \begin{aligned}
    &\lim_{\theta\to\infty} -\frac{2b_n+2a(\theta)^2}{\theta} = \frac{N-1+2\Delta}{2}, \\
    &\lim_{\theta\to\infty} \frac{4 u_n}{\theta^2} = \frac{n(N+1-n)(2n-N-1-2\Delta)(2n-N-1+2\Delta)}{4(2n-N)(2n-N-2)},
  \end{aligned}
  \end{align}
  and using instead \eqref{RecCoef2} for $N$ even, one obtains the recurrence coefficients 
  \begin{align}
  \begin{aligned}
    &\lim_{\theta\to\infty} -\frac{2b_n+2a(\theta)^2}{\theta} = \frac{N-1+2\Delta}{2}+\frac{(2\Delta-1)(N+1)}{4}\left(\frac{1}{2n-N-1}-\frac{1}{2n-N+1}\right), \\
    &\lim_{\theta\to\infty} \frac{4 u_n}{\theta^2} = \frac{n(N+1-n)(2n-N-2\Delta)(2n-N-2+2\Delta)}{4(2n-N-1)^2}.
  \end{aligned}
  \end{align}
  The weights of the para-Krawtchouk polynomials are similarily obtained from those of the para-Racah polynomials. Substituting \eqref{reparac} in \eqref{persymmetricweights} and \eqref{evenw} and taking the limit $\theta\to\infty$ gives 
  \begin{align}
  \begin{aligned}
  &w_{2s}  = \frac{(-\Delta-j)_N(-j)_s(-\Delta-j)_s}{2(-1)^{j+1}\binom{2j}{j}s!j!(\Delta)_{j+1}(1-\Delta)_s}, \quad
  &w_{2s+1}= \frac{(-\Delta-j)_N(-j)_s(\Delta-j)_s}{2(-1)^{j}\binom{2j}{j}s!j!(-\Delta)_{j+1}(1+\Delta)_s}, 
  \end{aligned}
  \end{align}
  for $N=2j+1$ and 
  \begin{align}
  \begin{aligned}
  &w_{2s}  = \frac{(-\Delta-j+1)_N(-j)_s(-\Delta-j+1)_s}{(-1)^{j}\binom{2j}{j}s!j!(\Delta)_{j}(1-\Delta)_s}, \quad
  &w_{2s+1}= \frac{(-\Delta-j+1)_N(-j+1)_s(\Delta-j)_s}{(-1)^{j+1}\binom{2j}{j}s!(j-1)!(-\Delta)_{j+1}(1+\Delta)_s}, 
  \end{aligned}
  \end{align}
  for $N=2j$. This corresponds to the features of the para-Krawtchouk polynomials defined in \cite{2012_Vinet&Zhedanov_JPhysA_45_265304} with parameter $2\Delta$. We have thus given a limiting relation between the para-Racah polynomials and the para-Krawtchouk polynomials.  

  \subsection{Connection with the dual-Hahn polynomials}
  If one sets $c = a + 1/2$, it can be seen from \eqref{bi-lattice} and \eqref{bi-lattice2} that the orthogonality points now form a single quadratic lattice of the form
  \begin{align}
  x_s = -\left(\frac{s}{2}+a\right)^2, \qquad s = 0,1,\dots,N. 
  \end{align}
  Furthermore, it can be shown that the para-Racah polynomials with $c = a + 1/2$ and $\alpha = 1/2$ reduce to dual-Hahn polynomials. It is straightforward to verify from the recurrence relation that
  \begin{align}
   \left(-\frac{1}{4}\right)^n P_n\left(-\frac{1}{4}y-a^2;N,a,a+\frac{1}{2},\frac{1}{2}\right) = r_n\left(y;N,\frac{4a-1}{2},\frac{4a-1}{2}\right)
  \end{align}
  where $r_n\left(y;N,\gamma,\delta\right)$ denotes monic dual-Hahn polynomials ­\cite{2010_Koekoek_&Lesky&Swarttouw}. 

  \subsection{Cases $\alpha = 0$ and $\alpha=1$}
  The para-Racah polynomials have special properties for $\alpha =0$ and $\alpha =1$. When $N=2j+1$, the recurrence coefficient $u_{j+1}$ becomes zero for both special values of $\alpha$. In view of \eqref{Jmatrix}, we can see that the Jacobi matrix $J$ now splits into two $(j+1)\times (j+1)$ tri-diagonal blocks: 
    \begin{align} \label{Jspec1}
  J=
  \begin{pmatrix}
  b_0& \sqrt{u}_1 &0 &0 & & & &  \\
  \sqrt{u}_1 & \ddots & \ddots &0 &&&& \\
   0 & \ddots & \ddots & \sqrt{u_j} &&&& \\
   0 &0 & \sqrt{u_j} & b_j &&&& \\1
   &&&& b_{j+1} & \sqrt{u_{j+2}} &0 &0 \\
   &&&& \sqrt{u_{j+2}} & \ddots & \ddots &0 \\
   &&&& 0 & \ddots & \ddots & \sqrt{u_N} \\
   &&&& 0 &  0     & \sqrt{u_N} & b_N 
  \end{pmatrix}.
  \end{align} 
  This means that the recurrence relation for the para-Racah polynomials splits into two and defines independent sets of orthogonal polynomials. The first set has polynomials from degree $0$ to $j$ and the second set has polynomials from degree $j+1$ to $N$. Furthermore, it can be seen from \eqref{odd_weights} that upon setting $\alpha =0$ or $\alpha=1$, half of the weights vanishes and the orthogonality grid thus reduce to a single quadratic lattice.
  
  When $N=2j$, similar properties arise. For $\alpha=0$, we have $u_{j+1}=0$ and the Jacobi matrix takes the form \eqref{Jspec1} with two blocks of dimension $(j+1)\times (j+1)$ and $(j)\times (j)$ respectively. Similarily, for $\alpha =1$, the Jacobi matrix becomes
  \begin{align} \label{Jspec2}
  J=
  \begin{pmatrix}
  b_0& \sqrt{u}_1 &0 &0 & & & &  \\
  \sqrt{u}_1 & \ddots & \ddots &0 &&&& \\
   0 & \ddots & \ddots & \sqrt{u_{j-1}} &&&& \\
   0 &0 & \sqrt{u_{j-1}} & b_{j-1} &&&& \\1
   &&&& b_{j} & \sqrt{u_{j+1}} &0 &0 \\
   &&&& \sqrt{u_{j+1}} & \ddots & \ddots &0 \\
   &&&& 0 & \ddots & \ddots & \sqrt{u_N} \\
   &&&& 0 &  0     & \sqrt{u_N} & b_N 
  \end{pmatrix}
  \end{align}
  where the first block has dimension $(j)\times (j)$ and the second block has dimension $(j+1)\times (j+1)$. Again, this means that the para-Racah polynomials splits into two sets of mutually orthogonal polynomials on a single quadratic bi-lattice. 

\section{Conclusion}

  To sum up, the para-Racah polynomials have been introduced and characterized. These are polynomials of a discrete variable that are orthogonal on finite quadratic bi-lattices. They are obtained from a novel truncation of the Wilson polynomials. The explicit expression, orthogonality property, recurrence relation and difference equation have been provided. The cases of odd and even numbers of cardinalities of the finite sets of polynomials must be distinguished. The para-Racah polynomials have the dual-Hahn polynomials as a special case and the para-Krawtchouk polynomials as a limiting case. (The latter are orthogonal on linear bi-lattices). Looking forward, we plan to obtain the $q$-generalization of the para-Racah polynomials by following an approach similar to the one adopted in this paper.

  In closing, we would like to point out that the para-Racah polynomials have already found applications in the general framework of quantum information \cite{2015_Lemay&Vinet&Zhedanov_PhysRevA}. Indeed, they have been used in the design of a spin chain that can perform the perfect transfer of quantum states and generate maximally entangled states. We trust that the para-Racah polynomials will prove to have many other uses.

\section*{Acknowledgments}
The authors would like to thank S. Tsujimoto for stimulating discussions. J.M.L. holds a scholarship from the Fonds de recherche du Qu\'ebec -- Nature et technologies (FRQNT). The research of L.V. is supported in part by NSERC. A. Z. wishes to thank the Centre de recherches math\'ematiques (CRM) for its hospitality.

\section*{References}
\bibliographystyle{elsarticle-num}
\bibliography{Bibliography_ParaRacah.bib}

\end{document}